%% file: mainv5.tex
\documentclass[reprint,aps,pre,floatfix]{revtex4-1}

\usepackage{xcolor}
\usepackage[normalem]{ulem}
\usepackage{amsmath,amsfonts,amscd,amssymb,graphicx}
\input{z_vardef}

\usepackage{lineno}
\usepackage{dcolumn}
\usepackage{overpic}
\usepackage{color}
\usepackage{graphicx}
\usepackage{mathrsfs}
\usepackage{subcaption}
\usepackage{adjustbox}
\usepackage{rotating}
\usepackage{float}
\makeatletter
\let\newfloat\newfloat@ltx
\makeatother
\usepackage{algorithm}
\usepackage{algpseudocode}
\algnewcommand\algorithmicinput{\textbf{Input:}}
\algnewcommand\Input{\item[\algorithmicinput]}

\begin{document}

\title{Sensing with shallow recurrent decoder networks}

\author{Jan P. Williams$^*$, Olivia Zahn$^{**}$ and J. Nathan Kutz$^\dag$}
 \affiliation{$^*$Department of Mechanical Engineering, University of Washington, Seattle, WA} 
 \affiliation{$^{**}$ Department of Physics, University of Washington, Seattle, WA}
 \affiliation{$^\dag$ Departments of Applied Mathematics and Electrical and Computer Engineering, University of Washington, Seattle, WA}

\begin{abstract}
Sensing is a universal task in science and engineering. Downstream tasks from sensing include inferring full state estimates of a system (system identification), control decisions, and forecasting. We propose a {\em SHallow REcurrent Decoder} (SHRED) neural network structure for sensing which incorporates (i) a recurrent neural network to learn a latent representation of the temporal dynamics of the sensors, and (ii) a shallow decoder that learns a mapping between this latent representation and the high-dimensional state space. SHRED enables accurate reconstructions with far fewer sensors, outperforms existing techniques when more measurements are available, and is more robust to random sensor placements. In the example cases explored, complex spatio-temporal dynamics are characterized with exceedingly limited sensors that can be randomly placed with minimal loss of performance.
\end{abstract}

\maketitle


\section{Introduction} \label{intro}

Emerging sensor technologies are transforming every science and engineering domain, with the quantity and quality of data collected also driving fundamental advances through data-science and machine learning methods.  In many areas of interest, measurements of the full-state are at best impractical and often impossible.  Thus sensors are commonly used to infer the current and future behavior of high-dimensional systems with a limited number of sensor locations.  With severely limited and noisy sensor measurements, this task is exceptionally difficult and frequently requires principled sensor placement schemes to yield faithful reconstructions \cite{manohar_data-driven_2018}. 
Accurate and robust reconstruction techniques are vital in enabling downstream tasks such as system identification, forecasting, and control. Existing data-driven techniques to learn mappings from sensor measurements to the full-state fall into two categories: those that rely only on static sensor measurements~\cite{erichson_shallow_2020, yu_flowfield_2019, callaham_robust_2019, bolton_applications_2019} and those that involve delay embedding techniques~\cite{takens_detecting_1981, gilpin_deep_2020, uribarri_dynamical_2022, young_deep_2023}. The former have been used extensively in the reconstruction of large spatio-temporal fields, while the latter have primarily been applied with the aim of reconstructing chaotic attractors. In this work, we apply such delay embedding techniques to high-dimensional spatio-temporal flows. The aim is not to replicate an attractor, but rather infer the full-state from sparse sensor measurements. As will be shown, not only is there a significant performance increase  in comparison to static sensing, but a minimal number of sensors, randomly placed, can be used.

\begin{figure*}[t]
    \centering
    \subfloat{%
      \begin{overpic}[width=0.96\linewidth]{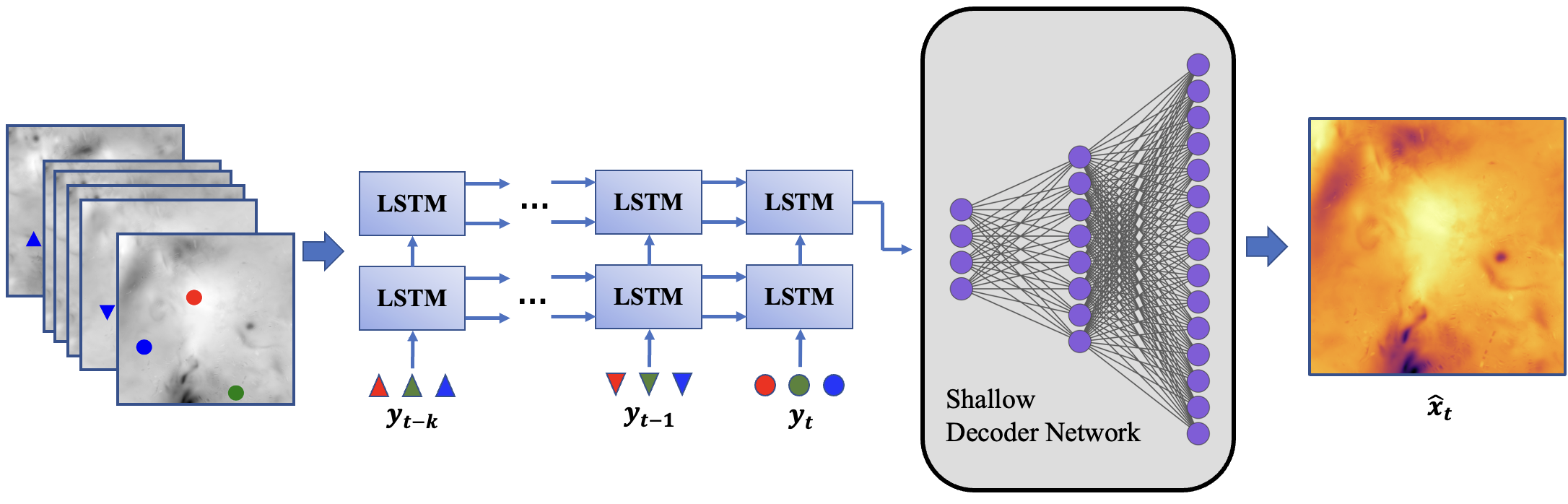}
        \put(2, 26){\textbf{A}}
      \end{overpic}
    }
    \vspace{0.2cm}
    \subfloat{%
      \begin{overpic}[width=0.32\textwidth]{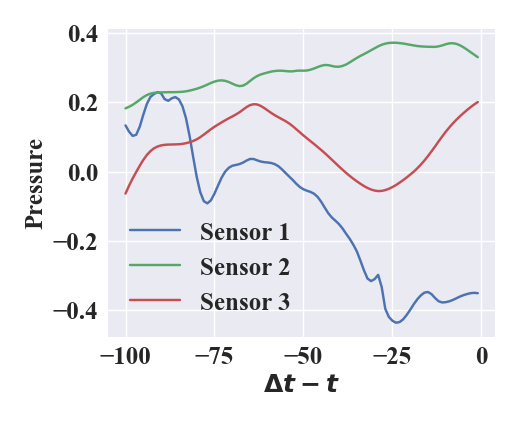}
        \put(3,75){\textbf{B}}
      \end{overpic}
    }
    \hfill
    \subfloat{
        \begin{overpic}[width=0.32\textwidth]{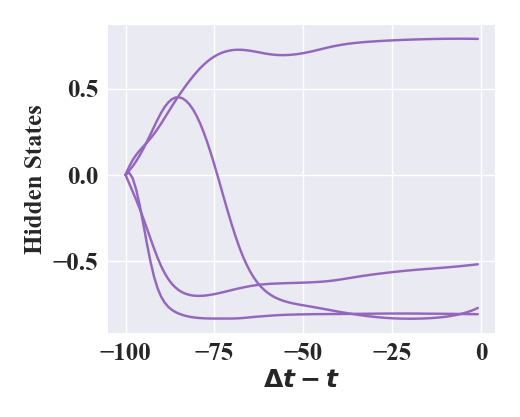}
        \put(3,75){\textbf{C}}
      \end{overpic}
    }
    \hfill
    \subfloat{
        \begin{overpic}[width=0.32\textwidth]{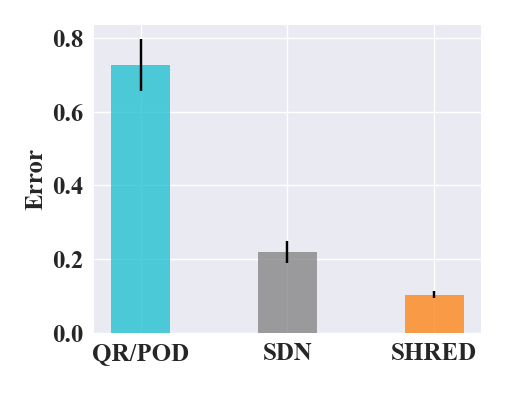}
        \put(3,75){\textbf{D}}
      \end{overpic}
    }
    \caption{Summary diagram of SHallow REcurrent Decoder networks (SHRED) for flow reconstruction from sensor measurements. \textit{A} A graphical representation of the SHRED method.  The mulitvariate time-series of sensor measurements, $\{ \by _i\}_{t-k}^t$, is fed into a stacked long short-term memory layer. The final output of the recurrent layer, $\mathbf{h}_t$, serves as an input to a fully-connected shallow decoder mapping from the hidden state to the high-dimensional field.  \textit{B} The time series of sensor measurements fed into the SHRED model. \textit{C} The evolution of the output hidden state generated by the input sequence of sensor measurements.  \textit{D} Reconstruction errors of traditional methods (QR/POD), shallow decoders (SDN), and SHRED on a turbulent flow when three sensors are 
    available.}
        \label{fig:summary-diagram}
\end{figure*}

The reconstruction of spatio-temporal dynamics from limited sensors relies on low-rank features of the data.  Computing low-rank embeddings of such high-dimensional data is often achieved through the {\em singular value decomposition} (SVD), also known as {\em proper orthogonal decomposition} (POD) \cite{kutz_data-driven_2013, brunton_data-driven_2019}.  The coefficients of the dominant correlated modes are determined by solving a linear inverse problem and the modes themselves serve as a linear map between measurements and the spatio-temporal state space.  Many of these linear methods are built upon the mathematical framework of gappy POD and have been successful across many disciplines \cite{everson_karhunenloeve_1995, drmac_new_2016, barrault_empirical_2004, chaturantabut_nonlinear_2010, manohar_data-driven_2018}.  More recently, shallow decoder networks (SDN) have leveraged advances in machine learning and AI to learn end-to-end, nonlinear maps between measurements and high-dimensional state spaces.  SDNs have been demonstrated to outperform their linear counterparts, particularly when the number of available sensors is exceedingly low \cite{erichson_shallow_2020, carter_data-driven_2021, sahba_wavefront_2022, williams_data-driven_2022}.  Both derivatives of gappy POD and SDNs rely on measurements at a single snapshot in time to reconstruct the corresponding high-dimensional state at that time.  

As a result of this static sensing, linear reconstruction methods are heavily reliant on optimally placed sensors for the inverse problem to be well-conditioned \cite{manohar_data-driven_2018, drmac_new_2016}.  In general, determining optimal sensor locations is a combinatorially hard problem and is infeasible for large search spaces.  For smaller spaces, there exist a number of well-known techniques for sensor placement \cite{boyd_convex_2004, joshi_sensor_2009, caselton_optimal_1984, krause_near-optimal_2008, lindley_measure_1956, sebastiani_maximum_2000, paninski_asymptotic_2005}, but they become computationally intractable in high-dimensional systems.  Instead, convex relaxations of the problem become necessary for obtaining approximate solutions. A greedy search using the QR decomposition with column pivoting offers one approximate solution and such methods have been critical in enabling accurate reconstructions of high-dimensional systems \cite{manohar_data-driven_2018, drmac_new_2016}.  However, greedy algorithms suffer from the fact that there is no guarantee that the sensor locations found are physically implementable; for instance, measuring sea-surface temperature in the middle of the Pacific Ocean is significantly more challenging than taking measurements along the coast.  Variations of greedy search algorithms can incorporate additional constraints, such as cost and sensor fidelity~\cite{clark_sensor_2020, clark_multi-fidelity_2021, clark_greedy_2019, manohar_optimized_2019}, but often at the expense of reconstruction accuracy.  The use of SDNs helps mitigate the dependence on sensor location, but still is greatly aided by principled sensor placement schemes \cite{williams_data-driven_2022}. 

The data-driven method presented here incorporates temporal  trajectories of sensor measurements to improve reconstruction accuracy, robustness to noise, and eliminate the need for optimally placed sensors. We use a type of recurrent network layer, {\em  long short-term memory networks} (LSTM) \cite{hochreiter_long_1997}, to process a time-series of sensor measurements.  The latent representation of the LSTM is the input into a fully-connected, shallow decoder for reconstruction.  We demonstrate that these {\em SHallow REcurrent Decoders} (SHRED) outperform existing linear and nonlinear techniques on three example datasets: a forced isotropic turbulent flow \cite{li_public_2008}, weekly sea-surface temperature \cite{reynolds_improved_2002}, and atmospheric ozone concentration \cite{bey_global_2001}.  In all cases, SHRED with as few as three randomly placed sensors achieves superior reconstruction accuracy than existing techniques using a far greater number of optimally placed sensors.  Moreover, because the inputs to a trained SHRED model consist of just a few sensor measurements, SHRED offers an {\em on-the-fly} compression for modeling physical and engineering systems.  Fig. \ref{fig:summary-diagram} gives a graphical summary of SHRED. The use of both constituents of SHRED in this paper, LSTMs and feed-forward shallow decoders, follows a long line of work demonstrating the efficacy of each component and their use together replicating chaotic attractors from partial observations \cite{uribarri_dynamical_2022, erichson_shallow_2020, kasai_deep_2021}. That being said, there exist a litany of different architectures that could be used to encode sensor trajectories (basic RNNs, GRUs \cite{chung_empirical_2014}, ESNs \cite{lukosevicius_practical_2012}, transformers \cite{vaswani_attention_2023}, etc.) or decode from a learned latent state (CNNs \cite{lecun_deep_2015}, transformers \cite{vaswani_attention_2023}, residual layers \cite{he_deep_2016}, etc.); the fundamental objective of SHRED is to combine implicit time-delay embeddings of sensor measurements with an expressive decoder for high-dimensional state reconstruction. The investigation of the many possible architectures to perform this task is an important area of future exploration. In this work, we have shown explicitly that the combination of the two key components, \textit{recurrence} and \textit{decoding}, allow for a successful mathematical scheme for sensing.

\section{Background}

\subsection{Proper orthogonal decomposition and QR}

Traditional techniques for high-dimensional state reconstruction from limited sensor measurements typically rely on the SVD to compute a low-rank embedding of the spatio-temporal dynamics.  Given $N$ snapshots of an $m$ dimensional state, we construct a data matrix $\bX = [ \bx _1 \; \bx_2 \; \dots \bx _N]$ where $\bx _i \in \mathbb{R}^m$ is the $i$-th temporal snapshot of the state.  The SVD factors the data matrix into the product of an orthogonal matrix, $\bU \in \mathbb{R} ^ {m \times m}$, a diagonal matrix with decreasing entries, $\mathbf{\Sigma} \in \mathbb{R} ^ {m \times N},$  and another orthogonal matrix $\bV \in \mathbb{R} ^{N \times N}.$ 
An optimal rank $r$ approximation of $\bX$ can be directly obtained through the SVD by retaining only the first $r$ columns of $\bU$,

\begin{equation}
    \bX \approx \hat \bX = \bU _r \mathbf{\Sigma} _r \bV _r ^T.
\end{equation}

Thus, each temporal snapshot $\bx _i$ can be approximated by a linear combination of the first $r$ columns of $\bU$, otherwise known as the dominant POD modes.  

POD based methods estimate a high-dimensional state, $\bx$, drawn from the same distribution as that of the training data by solving for these $r$ coefficients using a subsampling of the state,

\begin{equation}
    \by = \bC \bx \approx \bC \bU _r \bb,
\end{equation}

where $\bC$ is a ``sensing'' matrix consisting of rows of the $m \times m$ identity matrix and $\bb \in \mathbb{R}^r$.  A full-state estimate from a set of sensor measurements, $\by$, can then be obtained by

\begin{equation}
    \bx \approx \hat \bx = \bU _r (\bC \bU_r) ^{-1} \by.
    \label{eq:reconstruction}
\end{equation}

In practice, it is critical to design $\bC$ such that the inversion  in  (\ref{eq:reconstruction}) is well-conditioned.  An indirect bound on the condition number of $\bC \bU _r$ can be found by determining $\bC$ where $| \det (\bC \bU _r) |$ is maximized. That is, we seek 

\begin{align}
    \mathbf{C}.  = \argmax_{\bC} |\det \mathbf{C} \mathbf{U} _r| = \argmax_{\bC} \prod _i|\lambda _i (\bC U_r)| \notag \\ = 
    \argmax_{\bC} \prod _i\sigma _i (\bC U_r).
    \label{eq:qr_optimization}
\end{align}

Unfortunately, for large systems this $NP$-hard, combinatorial search becomes computationally intractable.  The QR decomposition with column pivoting of $\bU _r ^T$ offers an approximate greedy solution.  The decomposition utilizes Householder transformations to find a column permutation matrix $\mathbf{P}^T$ such that 

\begin{equation}
    \bU _r ^T \mathbf{P}^T = \bQ \bR
\end{equation}

where $\bQ$ is orthogonal and $\bR$ is upper triangular with decreasing, positive entries on the diagonal.  Because at each iteration of the algorithm the diagonal entry of $\bR$ is selected to be as large as possible and $\bQ$ is orthogonal, we have a a greedy maximization of $\det (\bU_r ^T \mathbf{P}^T) = \det (\mathbf{P} \bU_r).$ Thus, setting $\bC = \mathbf{P}$ can improve the reconstructions obtained by POD based methods.  We refer to reconstructions obtained in this manner as QR/POD.

\subsection{Shallow decoder networks}

Although QR/POD methods have demonstrated success in a variety of fields, the method can, at most, estimate the coefficients of the first $n$ POD modes if there are $n$ sensors available.  The included POD modes are insufficiently expressive to obtain accurate reconstructions when exceedingly few sensors are available.  This limitation motivates the development of more expressive, nonlinear methods for reconstruction.  Among these are shallow decoder networks (SDN) \cite{erichson_shallow_2020}.

As before, let

\begin{equation}
    \by = \bC \bx
\end{equation}

be a point sampling of a high-dimensional state $\bx.$  SDNs are shallow, fully-connected neural networks that learn a mapping from the space of sensor measurements, $\by$, back to $\bx.$  SDNs can be denoted as

\begin{equation}
    \mathcal{F}(\by; \mathbf{W}_{SD}) := R(\mathbf{W}^bR(\mathbf{W}^{b-1}\cdots R(\mathbf{W^1 s})))
    \label{eq:SDN}
\end{equation}

where $\by$ is the input sensor data, $R$ is a scalar, nonlinear activation function, $\mathbf{W} _{SD} = \{ \mathbf{W^i} \}_{i=1}^b$ are trainable weights, and $k$ denotes the number of layers in the SDN. Formally, we seek

\begin{equation}
    \mathcal{F} \in \argmin_{\widetilde{\mathcal{F}} \in  \mathscr{F}} \sum _{i=1}^N ||\bx_i - \widetilde{\mathcal{F}}(\by_i)||_2,
\end{equation}

so that reconstruction error is minimized over a set of $N$ training states $\{ \bx _i \}_{i=1}^N$.  The parameters of the network are randomly initialized and then trained by a gradient descent method.  Once the network is trained, reconstructions given subsequent sensor measurements are found by 

\begin{equation}
    \hat{\bx} = \mathcal{F} (\by).
\end{equation}

SDNs can be trained for any set of sensor measurements, and unlike QR/POD are not entirely reliant on principled sensor placement schemes.  Still, previous empirical work suggests that even neural-network based reconstructions can be improved through the use of QR placement \cite{williams_data-driven_2022}. 
Both SDNs and QR/POD have been used successfully in the reconstruction of high-dimensional dynamical systems, but also for broader classes of images~\cite{manohar_data-driven_2018}.  This extension is possible because both methods rely only on static, point measurements and each reconstruction is performed individually.  Dynamical systems, on the other hand, inherently depend on the temporal evolution of the state, and incorporating this dependence into a method for reconstruction offers a natural improvement upon existing techniques.  Our work includes these dynamics through the use of trajectories of sensor measurements.

\section{Shallow recurrent decoder networks}

\subsection{SHRED and separation of variables}

The SHRED architecture is based upon the separation of variables technique for solving linear partial differential equations (PDEs)~\cite{folland1995introduction}.  Separation of variables assumes that a solution can be separated into a product of spatial and temporal functions $u(x,t)= T(t) X(x)$.  This solution form is then used to reduce the PDE into ordinary differential equations:  one for time $T(t)$ and one for space $X(x)$.  Such a decomposition also constitutes the underpinnings of spectral methods for the numerical solution of the PDE, linear or nonlinear~\cite{kutz_data-driven_2013}.  
Consider the constant coefficient linear PDE 
\begin{equation}
 \dot{u} = {\cal L} (\partial_x, \partial^2_x, \cdots ) {u}
 {\label{eq:linearPDE}}
\end{equation}
where $u(x,t)$ specifies the spatio-temporal field of interest. Typically the initial condition (IC) and boundary conditions (BCs) are given by
\begin{subequations}
\begin{align}
\text{IC:}&&\quad  u(x,0)=u_0(x) \\
\text{BCs:}&&\quad  \alpha_1 u (0,t) +\beta_1 u_x(0,t) = g_1(t) \,\, \text{and} \nonumber \\  && \alpha_2 u(L,t) + \beta_2 u_x(L,t) = g_2(t).
\end{align}
\label{eq:ICBC}
\end{subequations}
This may be generalized to systems of several spatial variables, or a system with no time dependence.  The linear operator ${\cal L}$ specifies the spatial derivatives, which in turn model the underlying physics of the system.  Simple examples of ${\cal L}$ include ${\cal L}= c \partial_x$ (the one-way wave equation) and ${\cal L}= \kappa \partial^2_x$ (the heat equation)~\cite{kutz_data-driven_2013}.  

The earliest solutions of linear PDEs assumed separation of variables whereby $u(x,t)=\exp(\lambda t) X(x)$ was a product of a temporal (exponential) function multiplied by a spatial function.  The parameter $\lambda$ is in general complex.  This gives the eigenfunction solution of (\ref{eq:linearPDE}) to be
\begin{equation}
    u(x,t) = \sum_{n=1}^{N} a_n \exp(\lambda_n t) \phi_n(x)
    \label{eq:ef}
\end{equation}
where $\phi_n(x)$ are the eigenfunctions of the linear operator and $\lambda_n$ are its eigenvalues (${\cal L}\phi_n(x) = \lambda_n \phi_n(x)$).  Here a finite dimensional approximation $N$ is assumed, which is standard in practice for numerical evaluation.  

Typically, initial conditions $u(x,0)=u_0(x)$ are imposed in order to uniquely determine the coefficients $b_n$.  Specifically, at time $t=0$ (\ref{eq:ef}) becomes
\begin{equation}
    u_0(x) = \sum_{n=1}^{N} a_n  \phi_n(x) .
\end{equation}
Taking the inner product of both sides with respect to $\phi_m(x)$ and making use of orthogonality gives
\begin{equation}
   a_n =\langle u_0(x), \phi_n(x) \rangle
\end{equation}
SHRED instead has measurements at a single spatial (sensor) location $x_s$, but with a temporal history.  Thus if SHRED has, for example, $N$ temporal trajectory points, this gives at each time point of the measurement:
\begin{equation}
    u(x_s, t_j) = \sum_{n=1}^{N} a_n 
    \exp(\lambda_n t_j) \phi_n(x_s) 
    \,\,\,\,\,\, \mbox{for} \,\,\, j=1,2,\cdots N .
\end{equation}
This results in $N$ equations for the $N$ unknowns $a_n$. Specifically, the $N\times N$ system of equations ${\bf A} {\bf x} = {\bf b}$ is prescribed by the vector components $x_k=a_k$ and $b_k = u(x_s,t_k)$ and matrix components $(a_{kj}) = \exp(\lambda_k t_j) \phi_k(x_s)$.  As with the initial condition (\ref{eq:ICBC}a), the time trajectory of measurements at a single location uniquely prescribes the solution. This analysis can be generalized to allow for spatially mobile sensors or include multiple sensor measurements at a single time point. The former case is addressed in \cite{ebers2023leveraging} and generally requires that the sensor’s movement be periodic. In the latter case, if there are two measurements at a given time $t_j$, then only $N/2$ trajectory points are needed to uniquely determine the solution. Likewise, three sensor measurements at a given time requires only $N/3$ trajectory points. As a caveat, it should be noted that sensors must actually measure the physics in order to perform a reconstruction, i.e. if a shock wave has not hit the sensor yet, one cannot construct the spatio-temporal dynamics.  Moreover, SHRED cannot reconstruct regions of space where the dynamics are statistically independent.

Thus, temporal trajectory information at a single spatial location, or with a moving sensor, can be equivalent to knowing the full spatial field at a given point in time. The motivation of SHRED is a nonlinear generalization of this concept. The proposed architecture has the form $u = X(T(\{y_i\}^{i=t}_{i=t-k}))$. In separation of variables, one assumes the form $u = X(x)T(t)$ which allows for the construction of differential equations for time and space separately:  $dT/dt = c T$ and ${\cal L}X=cX$.  However, they are not independent, but are related through the constant $c$.  Specifically, the solution $T(t)$ is constrained to constants that are the eigenvalues of $X(x)$. In a SHRED model, an analogous result is obtained through the joint training of $X$ and $T$. Of course, as a generalization to nonlinear dynamics, rigorous theoretical bounds of SHRED are difficult to achieve, much like analytic and numerical solutions are difficult to rigorously bound in computational PDE settings.

That being said, standard numerical techniques for non-linear PDEs can also be obtained with time-series of sparse sensor measurements. Consider a spectral solution technique \cite{kutz_data-driven_2013} whereby numerical solutions are approximated by a spectral basis 
 \begin{equation}
     u(x,t) = \sum_{n=1}^{N} a_n (t) \phi_n(x) .
     \label{eq:spectral}
\end{equation}
This spectral decomposition turns the PDE into a system of $N$ coupled ordinary differential equations with constants of integration that can be uniquely determined by requiring the solution satisfy $N$ temporal trajectory points.

\subsubsection*{Coupled PDEs}

We can also consider coupled, constant coefficient linear PDEs.  For example, the coupled system 
\begin{subequations}
    \begin{eqnarray}
 \dot{u} = {\cal L}_1  {u}
 + {\cal L}_2  {v} {\label{eq:linearPDE2a}}\\
 \dot{v} = {\cal L}_3  {u}
 + {\cal L}_4 {v}
 {\label{eq:linearPDE2b}}
\end{eqnarray}
\end{subequations}
where $u(x,t)$ and $v(x,t)$ specifies the spatio-temporal fields of interest. The PDEs can be instead be written in the form
\begin{equation}
    \ddot{u} = {\cal L}_1 \dot{u}
    +{\cal L}_2 {\cal L}_3 u
    +{\cal L}_2 {\cal L}_4 \left(
     {\cal L}_2^{-1} (\dot{u} - {\cal L}_1 {u})
    \right)
\end{equation}
where (\ref{eq:linearPDE2a}) is differentiated with respect to time and (\ref{eq:linearPDE2b}) is used in order to write the PDEs as a function of $u(x,t)$ alone.  Thus knowledge of the field $u(x,t)$ alone is capable of constructing the solution fields $u(x,t)$ and $v(x,t)$.  For this second-order (in time) PDE, both an initial condition $u(x,0)$ and an initial {\em velocity} require specification 
$\dot{u}(x,0)$ in order to uniquely determine the solution.  As with the previous arguments, a time trajectory embedding of $2N$ measurements can be used to uniquely determine the solution.
\subsection{SHRED architecture}
To model temporal dynamics, we rely on recurrent neural networks (RNN) \cite{rumelhart_learning_1987} and, more specifically, long short-term memory networks (LSTM).  The most basic RNN accepts as an input a sequence of vectors $\{ \bv_i \} _{i=1}^t$ and outputs a single hidden state, $\mathbf{h}_t.$  This can be expressed through the recursion relation

\begin{equation}
    \mathbf{h}_t = R(\bW _c \bv _t + \bW _r \mathbf{h}_{t-1} + \bb _r)
\end{equation}

where $R$ is a scalar, nonlinear activation function, $\bW _c$, $\bW _r$ and $\bb_r$ are trainable weights, and $\mathbf{h_0} = \begin{bmatrix}
0 & \cdots & 0
\end{bmatrix}^T$.  Generic RNNs are notorious for suffering from the vanishing gradient problem, causing them to have difficulty in identifying long term dependencies.  LSTMs address this issue through the introduction of a so-called ``gradient super-highway.''  Instead of outputting a single hidden state, LSTMs incorporate an additional ``cell state'' which only undergoes minor, pointwise operations, allowing gradients to flow easily from many time steps in the past.  LSTMs are perhaps the most commonly used class of RNN at the time of writing, due to their ability to learn long and short term dependencies \cite{lindemann_survey_2021}.  They are frequently used in time-series forecasting, natural language processing, and video analysis, among many other domains \cite{lindemann_survey_2021}.  For a more complete discussion of LSTMs, we direct the reader to \cite{hochreiter_long_1997}.  The recursion relation of an LSTM can be written as 

\begin{gather}
    \bh _t = \sigma \left( \bW _o \begin{bmatrix}  \bh_{t-1}, \bv_t \end{bmatrix} + \bb _o \right) \odot \tanh(\bc _t)\\
    \bc_t = \sigma \left( \bW _f \begin{bmatrix}
        \bh _{t-1}, \bv_t  
    \end{bmatrix} +\bb _f \right) \odot \bc _{t-1} \\ 
    + \sigma \left( \bW _i \begin{bmatrix}
        \bh _{t-1}, \bv_t
    \end{bmatrix} +\bb _f \right) \odot \tanh \left( \bW_g \begin{bmatrix}
        \bh_{t-1}, \bv_{t}
    \end{bmatrix} + \bb _g \right) \notag
\end{gather}

for trainable weights and biases $\bW _{RN} = \{\bW _o, \bW _f, \bW _i, \bW _g, \bb _o, \bb _f, \bb _i, \bb _g \}$.  Both $\sigma$, the sigmoid function, and $\tanh$ operate pointwise.  Note that $\mathbf{h}_t$ can be written as function in terms of only $\{ \bv _i \} _{i=1}^t$ and parameterized by the trainable weights and biases,

\begin{equation}
    \mathbf{h}_t = \mathcal{G} (\{ \bv _i \} _{i=1}^t; \bW_{RN}). 
    \label{eq:LSTMfunc}
\end{equation}

We now introduce the {\em SHallow REcurrent Decoder} (SHRED), which merges an LSTM with the SDN architecture.  Suppose we have a multivariate time-series of a high-dimensional state $\{ \bx_i \}_{i=1}^T$ and corresponding measurements of the system $\{ \by _i \}_{i=1}^T.$  Rather than reconstruct $\bx _i$ using only $\by _i$, as is done by QR/POD and SDN, we let an LSTM learn a latent representation using the previous $k$ sets of sensor measurements.  This latent representation is then used by a shallow decoder to reconstruct the high-dimensional state.  The SHRED architecture can be written as 

\begin{equation}
    \mathcal{H} \left( \{ \by _i \}_{i=t-k}^t \right) = \mathcal{F} \left( \mathcal{G} (\{ \by _i \}_{i=t-k}^t; \bW_{RN}); \bW_{SD}) \right)
\end{equation}

using eqs. \ref{eq:SDN} and \ref{eq:LSTMfunc}.  As in the case of SDNs, we use a gradient descent method to train the network weights to minimize reconstruction loss, 

\begin{equation}
    \mathcal{H} \in \argmin_{\widetilde{\mathcal{H}} \in  \mathscr{H}} \sum _{i=1}^N ||\bx_i - \widetilde{\mathcal{H}}\left( \{ \by _i \}_{i=t-k}^t \right)||_2.
\end{equation}

$\mathcal{H}$, so defined, represents a map from sensor trajectories to the full-state.

\section{Results}

\begin{figure*}
    \centering
    \begin{minipage}{0.52\linewidth}
        \subfloat{
      \begin{overpic}[width=\linewidth]{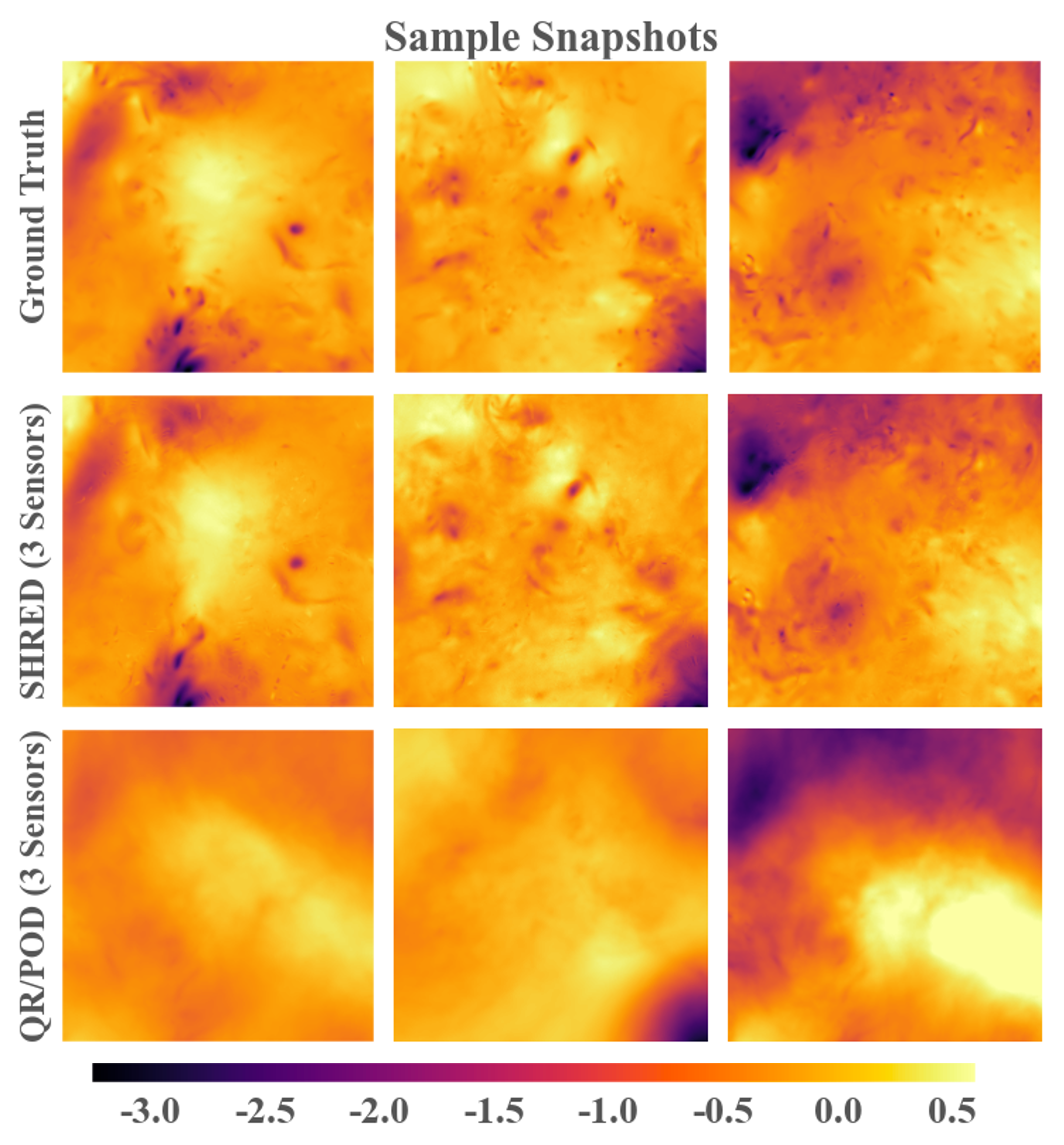}
        \put(2,95){\textbf{A}}
      \end{overpic}}
    \end{minipage}
    \begin{minipage}{0.415\linewidth}
    \subfloat{
    \begin{overpic}[width=0.98\linewidth]{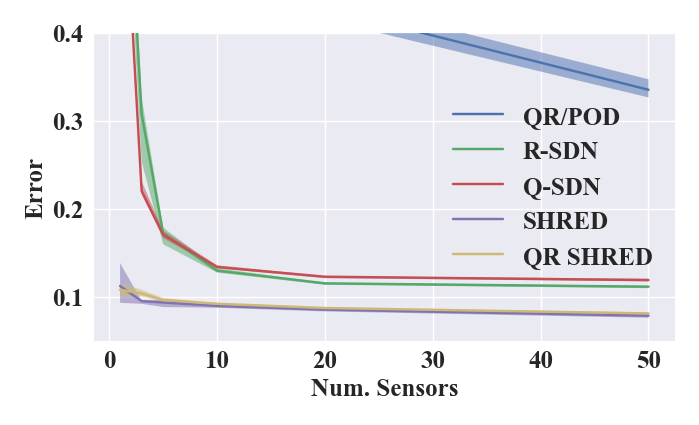}
      \put(2,55){\textbf{B}}
    \end{overpic}}
        
       \subfloat{
    \begin{overpic}[width=0.98\linewidth]{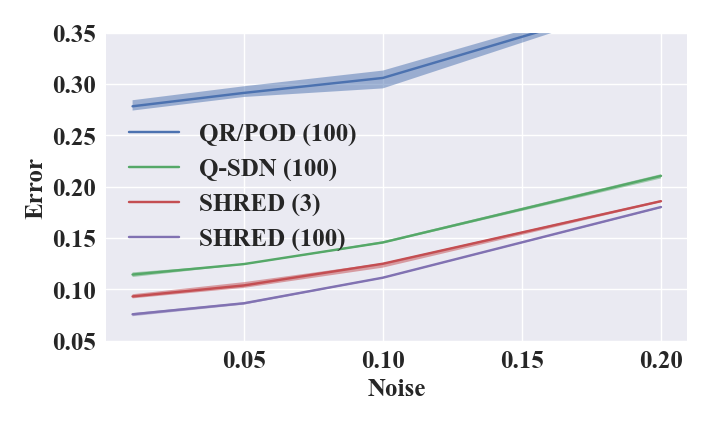}
      \put(2,55){\textbf{C}}
      \put(61,4.25){$(\alpha)$}
    \end{overpic}}
    \end{minipage}
    \caption{Turbulent flow summary performance. \textit{A} Example reconstructions obtained via SHRED and QR/POD of a turbulent flow when three sensors are available. Ground truth included for comparison.  \textit{B} Reconstruction errors of the current state-of-the-art static sensing methods and SHRED with a varying number of available sensors.  The solid lines denote the median error from 32 trained estimators and the shaded region denotes the interquartile range.  \textit{C} Performance of reconstruction methods in the presence of varying levels of added Gaussian white noise.  The added noise has mean zero and standard deviation equal to $\alpha$ times the mean absolute value of the field.  The number in parentheses denotes the number of available sensors.}
    \label{fig:ISOPanel}
\end{figure*}
\begin{figure*}[!ht]
    \centering
    \begin{minipage}{0.52\linewidth}
        \subfloat{
      \begin{overpic}[width=\linewidth]{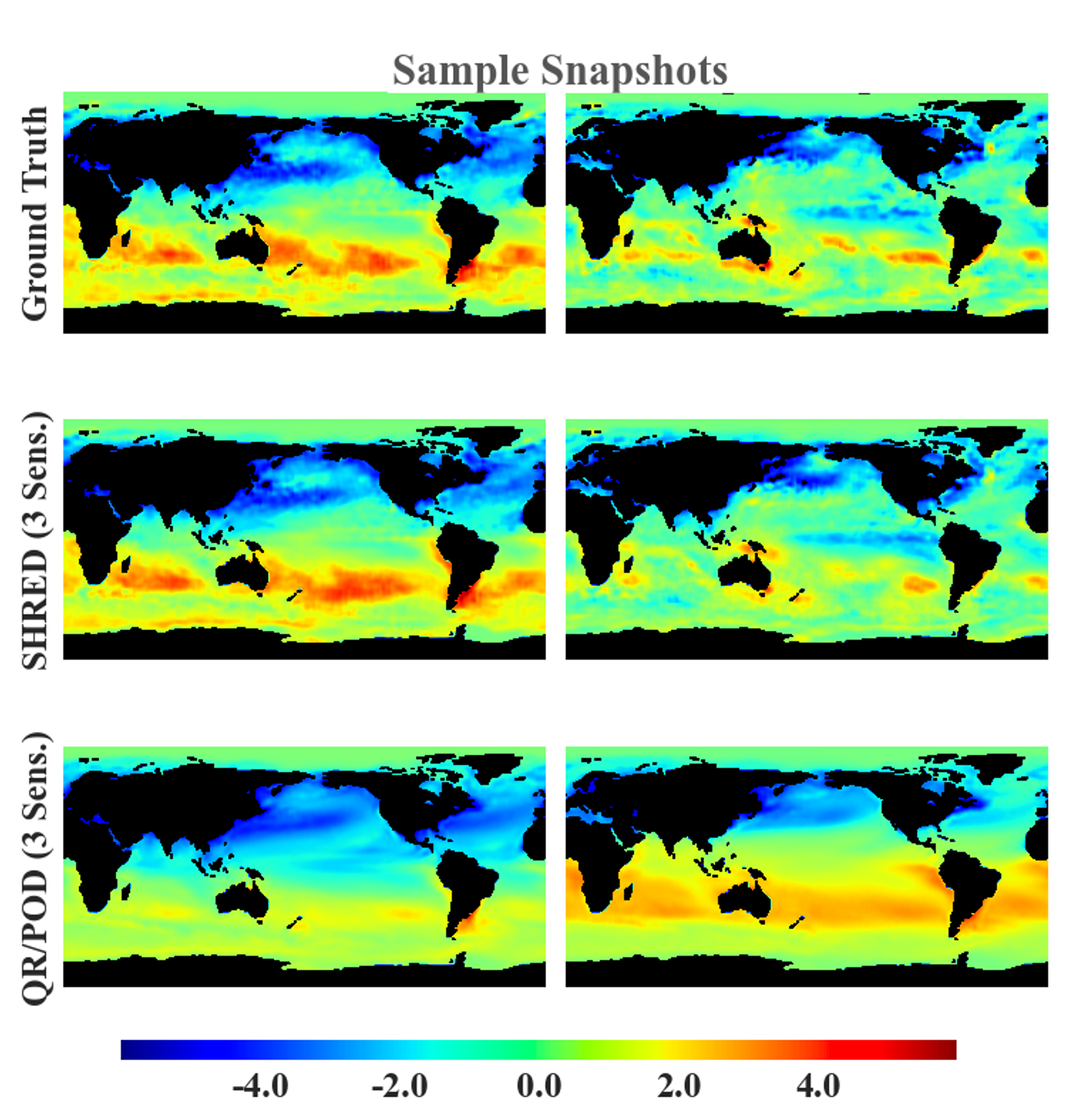}
        \put(2,95){\textbf{A}}
      \end{overpic}}
    \end{minipage}
    \begin{minipage}{0.415\linewidth}
    \subfloat{
    \begin{overpic}[width=0.98\linewidth]{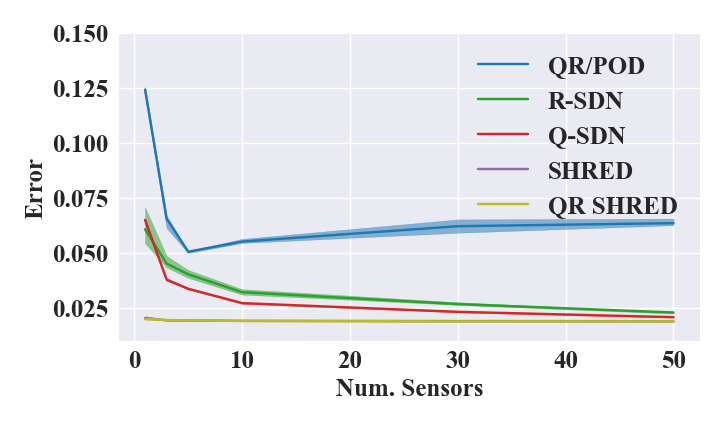}
      \put(2,55){\textbf{B}}
    \end{overpic}}
        
       \subfloat{
    \begin{overpic}[width=0.98\linewidth]{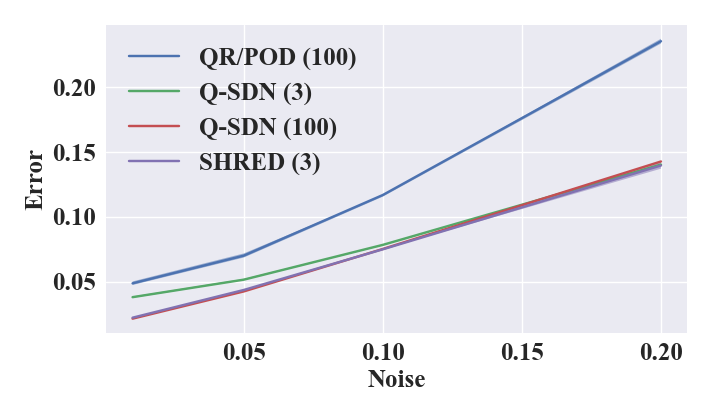}
      \put(2,55){\textbf{C}}
      \put(61,4.25){$(\alpha)$}
    \end{overpic}}
    \end{minipage}
    \caption{SST summary performance. \textit{A} Example reconstructions obtained via SHRED and QR/POD of sea-surface temperature when three sensors are available. Ground truth included for comparison.  \textit{B} Reconstruction errors of the current state-of-the-art static sensing methods and SHRED with a varying number of available sensors.  The solid lines denote the median error from 32 trained estimators and the shaded region denotes the interquartile range. \textit{C} Performance of reconstruction methods in the presence of varying levels of added Gaussian white noise.  The added noise has mean zero and standard deviation equal to $\alpha$ times the mean absolute value of the field.  The number in parentheses denotes the number of available sensors.}
    \label{fig:SSTPanel}
\end{figure*}
\begin{figure*}[t]
    \centering
    \begin{minipage}{0.52\linewidth}
        \subfloat{
      \begin{overpic}[width=\linewidth]{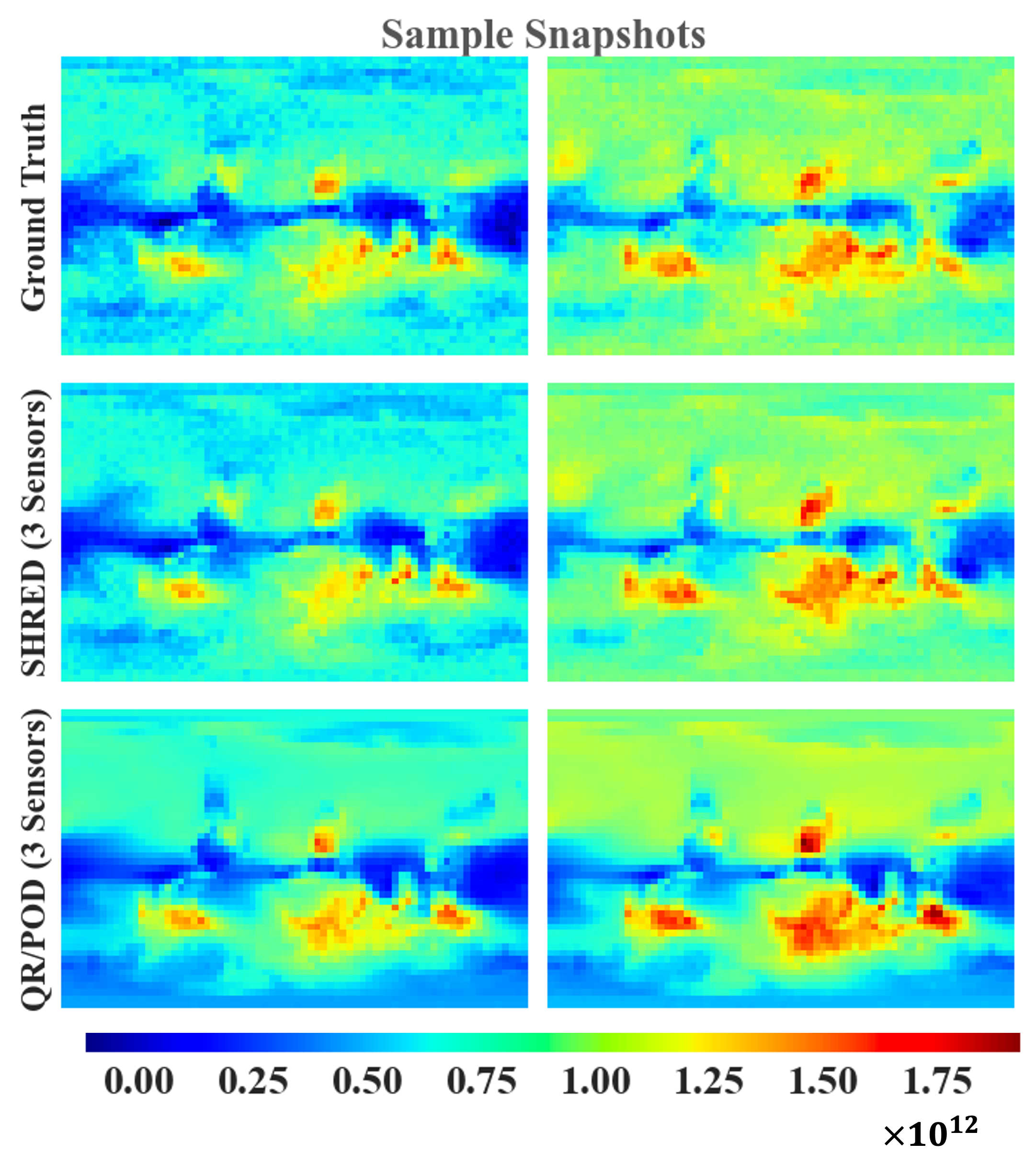}
        \put(2,95){\textbf{A}}
      \end{overpic}}
    \end{minipage}
    \begin{minipage}{0.415\linewidth}
    \subfloat{
    \begin{overpic}[width=0.98\linewidth]{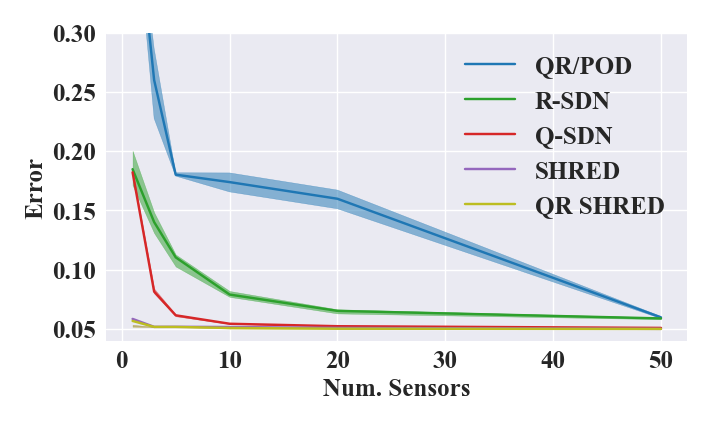}
      \put(2,55){\textbf{B}}
    \end{overpic}}
        
       \subfloat{
    \begin{overpic}[width=0.98\linewidth]{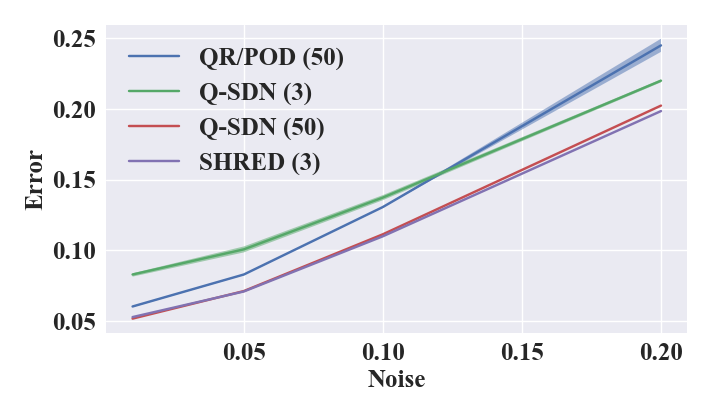}
      \put(2,55){\textbf{C}}
      \put(61,4.25){$(\alpha)$}
    \end{overpic}}
    \end{minipage}
    \caption{Ozone summary performance. \textit{A}  Example reconstructions obtained via SHRED and QR/POD of atmospheric ozone concentration when three sensors are available.  Only one of thirty elevations is shown.  Ground truth included for comparison.  \textit{B} Reconstruction errors of the current state-of-the-art static sensing methods and SHRED with a varying number of available sensors.  The solid lines denote the median error from 32 trained estimators and the shaded region denotes the interquartile range. \textit{C} Performance of reconstruction methods in the presence of varying levels of added Gaussian white noise.  The added noise has mean zero and standard deviation equal to $\alpha$ times the mean absolute value of the field.  The number in parentheses denotes the number of available sensors.}
    \label{fig:AO3Panel}
\end{figure*}

\subsection{Reconstruction of high-dimensional spatio-temporal fields}
\label{sec:results1}

We demonstrate the performance of SHRED trained with the Adam optimizer \cite{kingma_adam_2017} on three example datasets. In each case, we compare the reconstruction error obtained by SHRED with randomly placed sensors to SHRED with QR placed sensors (QR SHRED), shallow decoder networks with randomly placed sensors (R-SDN), shallow decoder networks with QR placed sensors (Q-SDN), and linear reconstructions with QR placed sensors (QR/POD). The considered reconstruction error is defined to be the normalized mean square error over each state in a test set \cite{erichson_shallow_2020},

\begin{equation}
    \text{Error} = \frac{|| \mathcal{H} \left( \{\by _j \}_{i-k}^i\right) - \bx_i ||_2}{|| \bx _i ||_2}.
\end{equation}

Because SHRED models rely on trajectories of sensor measurements to perform state estimation, we truncate each dataset to reconstruct only the final $N - k$ temporal snapshots, where $N$ is the initial number of samples and $k$ is the length of the utilized trajectories. This length can be viewed as a hyper-parameter that can be tuned according to the data. Detailed network parameters and training protocols can be found in the SI.  The used parameters of each model were selected via hyperparameter search. Finally, we note that in this section, training, validation, and test samples are temporally interspersed.  In later sections, we demonstrate results in the case that training and test sets are temporally distinct, although in all cases we assume no errors in sensor drift or synchronisation. The study of our method under these more difficult conditions is an important area of future work.

The first application we consider is the pressure field of a forced isotropic turbulent flow from the Johns Hopkins Turbulence Database \cite{li_public_2008}. The flow was generated by direct numerical simulation using $1024^3$ nodes and the pseudo-spectral method. We select a 350 by 350 two dimensional cutout of computed pressure over 1667 evenly spaced temporal snapshots from the three dimensional simulation.  The temporal spacing of the snapshots is approximately 0.006 seconds, while the DNS used to generate the data utilized time steps of 0.0002 seconds.  We seek to reconstruct these high-dimensional states from trajectories of point measurements of the states. There are many options for determining the optimal parameter $k$,  where $k\Delta t$ is the maximal time-delay utilized, for reconstructions generated by SHRED. In this case, a simple hyperparameter search over $k \in \{10, 25, 50, 75, 100\}$ was used with the number of sensors held constant at 3. From these possible $k$, the best performing value of $k = 100$ was selected.  Correspondingly, the final 1567 temporal snapshots of the pressure field are randomly split into training, validation, and test sets consisting of 1100, 234, and 233 snapshots, respectively. For each considered number of sensors, we generate 32 reconstructions with all considered methods. We plot the median performance and denote the interquartile range by the shaded region. These results are shown in panel B of Fig. \ref{fig:ISOPanel}.  Even with only a single, randomly placed sensor, the reconstructions obtained by SHRED yield significantly lower error than competing methods with as many as 50 sensors. Moreover, placement via the greedy QR algorithm appears to have a negligible impact on reconstructive performance. Panel A of Fig. \ref{fig:ISOPanel} shows sample reconstructions obtained by SHRED and QR/POD with three sensors. While QR/POD is only able to identify large scale features, SHRED accurately reconstructs fine grain features as well.  Numerically, with 1 sensor the median computed error of QR/POD is 0.89 while with SHRED the median computed error is 0.11, an approximately 8 times improvement in performance. The poor performance of QR/POD in this example highlights a well-known limitation of gappy-POD based techniques. The turbulent flow is both tri-periodic and nonstationary, presenting challenges that POD is fundamentally unable to handle. Such techniques remain widespread in part because of their relative interpretability and in part because there remains a lack of methods for inferring high-dimensional states from sparse sensor measurements alone. The comparison to POD illustrates the need for nonlinear methods in such cases.

In the vast majority of real world applications, and in contrast to numerical simulations, data is often corrupted by noise.  As a result, methods for state estimation must exhibit resilience to noise.  To measure this resilience, we corrupt the data with Gaussian noise of mean zero and standard deviation of $\alpha \times |\bar{x}|$, where $|\bar{x}|$ is the average absolute value over all points in all snapshots of the training set. We then generate 32 state estimates exactly as before using all considered methods. The resulting median reconstruction error and interquartile range for varying $\alpha$ is shown in panel C of Fig. \ref{fig:ISOPanel}.  Again, SHRED outperforms competing techniques using a far greater number of sensors.

The second application we examine is that of sea-surface temperature (SST).  We consider weekly mean sea-surface temperature from the years 1992 to 2019 as reported by NOAA  \cite{reynolds_improved_2002}.  Unlike the previous example of a simulated turbulent flow, SST is a sensor generated dataset for which governing equations are not known and thus represents a more practical application of SHRED.  The data consists of 1400 snapshots of a 180 by 360 grid, of which 44,219 spatial locations correspond to the sea-surface.  Each snapshot corresponds to the weekly mean sea-surface temperature for the respective week.  In this case, it is easy to make physical arguments for determining $k$, the number of lags used by SHRED. Since each data snapshot represents weekly mean sea-surface temperature, we expect a strong signal with period 52 corresponding to annual fluctuations in temperature. As a result, we select $k=52$ to capture this annual variation. We allow the input trajectories to SHRED to have a length of 52, corresponding to one year of measurements.  The final 1348 samples are divided into training, test, and validation sets consisting of 1000, 174, and 174 snapshots, respectively.  Analogous to Fig. \ref{fig:ISOPanel}, Fig. \ref{fig:SSTPanel} shows example reconstructions and reconstruction error distributions of SHRED and existing techniques both in the presence and absence of added Gaussian noise.  The contour plots shown in panel A of Figure 3 show the ground truth and reconstructions with the mean field subtracted. We choose to present this visualization because in comparison to other example datasets, the global mean-field dominates and the mean subtracted contour plots more clearly demonstrate the difference between QR/POD and SHRED. Again, the reconstructions obtained by SHRED are both visually and empirically superior when few sensors are available. When one sensor is available, the median error from QR/POD (0.12)  is approximately 6 times greater than that of SHRED (0.02).  Again, the reconstructions obtained by SHRED are both visually and empirically superior, while requiring far fewer sensors, which can be randomly placed without a loss of accuracy.


The last system for which we consider the performance of SHRED is a simulation of atmospheric chemistry.  Chemical transport models (CTM) simulate the evolution of an ensemble of interacting chemical species through a transport operator \cite{velegar_scalable_2019, bey_global_2001}

\begin{equation}
    \frac{\partial n_i}{\partial t} = - \nabla \cdot (n_i \bU)
\end{equation}

and a chemical operator

\begin{equation}
    \frac{dn_i}{dt} = (P_i -L_i) (\mathbf{n}) + E_i - D_i,
\end{equation}

where each entry $\mathbf{n} = \begin{bmatrix}
    n_1 & n_2 & \cdots n_K
\end{bmatrix}^T$ represents the number density of a specific chemical species, $\bU$ is the wind vector, $(P_i - L_i)(\mathbf{n})$ is the local chemical production and loss term, $E_i$ the emission rate of a species, and $D_i$ the deposition rate.  The output of a CTM, then, consists of concentrations for $K$ chemical species for a grid of latitudes, longitudes, and elevations.  The data we consider is drawn from the work of Velegar et al. \cite{velegar_scalable_2019} and contains simulated atmospheric ozone concentration generated by the CTM software GEOS-Chem \cite{bey_global_2001} over the course of a year with dynamical time steps of 20 minutes.  The accessed data from \cite{velegar_scalable_2019} is a compressed SVD representation using the first 50 POD modes.  The decompressed data matrix consists of 26,208 snapshots of a 46 by 72 by 30 (latitude, longitude, elevation) grid with temporal spacing of 20 minutes, from which we further downsample to obtain 2,600 global ozone concentration fields with temporal spacing of roughly 200 minutes.  To account for the fact that the data matrix is inherently rank 50, we add Gaussian noise with standard deviation $\alpha =0.05$ to the data before performing any analyses. Identical to the turbulent flow example, $k = 100$ for the reconstructions generated by SHRED was chosen by a hyperparameter search. The resulting snapshots are divided into training, validation, and test sets consisting of 1900, 300, and 300 entries, respectively.

Fig. 4 shows the performance of SHRED and existing techniques on this atmospheric ozone experiment. With 1 sensor, the median computed error of QR/POD (0.44) is roughly 7 times that of SHRED (0.06).  However, in this case when many sensors are available ($>50$) QR/POD performs comparably to nonlinear reconstructions.  This performance is an artifact of the use of a compressed, rank 50 representation of the data.

We note that for both the atmospheric ozone concentration and turbulent flow results, we have trained on synthetic data from validated high-fidelity models. As full-state measurements of physical high-dimensional systems are frequently impossible, such training procedures are common in data-driven reconstruction methods. However, the performance of any machine learning model trained using synthetic data is critically reliant upon the fidelity of the model used to generate the synthetic data. In real-world applications devoid of a well-validated numerical model or a training set of high-dimensional state measurements, it is unlikely that any ML model mapping from sparse sensor measurements to the full-state will perform well. One manner in which this simulation-to-reality gap can be addressed is by performing an expensive data collection period using a large number of deployed sensors. A shallow recurrent decoder could then be used to map from a small number of sensors to the larger set of sensor measurements as opposed to the full state. After the training period, costs can be reduced by removing many of the sensors that were initially deployed and interpolating between the sensor locations to which the map has been learned.

\begin{figure}[!ht]
    \centering
    \subfloat{
    \begin{overpic}[width=0.93\linewidth]{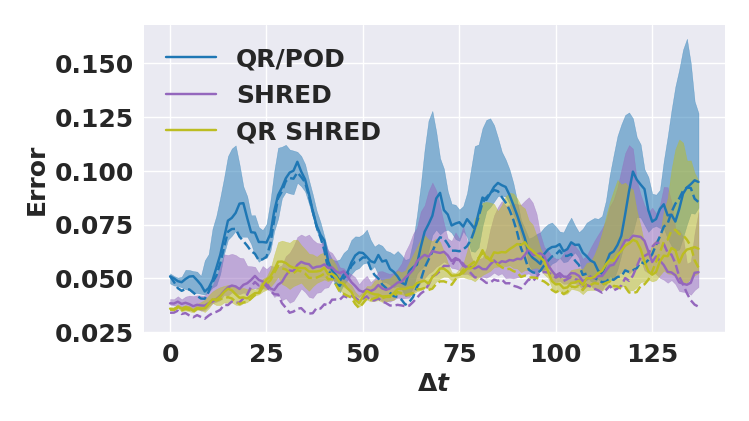}
      \put(2,55){\textbf{A}}
    \end{overpic}}
    
    \subfloat{
    \begin{overpic}[width=0.93\linewidth]{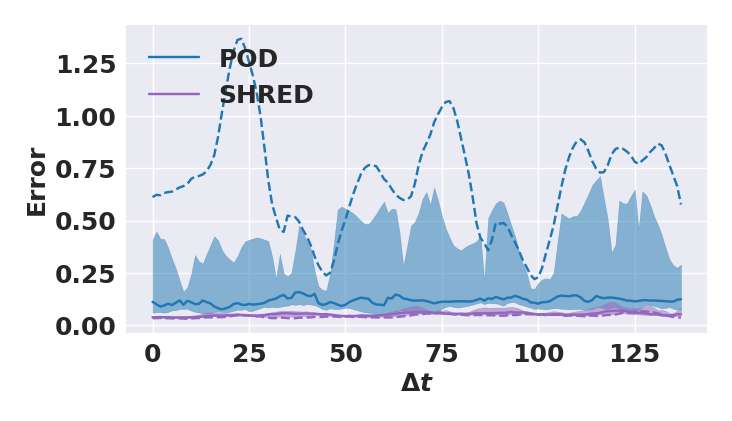}
      \put(2,55){\textbf{B}}
    \end{overpic}}
    
    \caption{Forecasting error for sea-surface temperature with $\Delta t$ as one week. A LSTM is trained to forecast QR placed, panel \textit{A}, or randomly placed, panel \textit{B}, sensor measurements which are subsequently used to perform reconstructions using SHRED and gappy POD.  16 forecasts are performed and the median error is denoted by the solid lines, with the 25th and 75th percentiles of reconstruction error defining the shaded region.  The dashed line represents an ensembled forecast.}
    \label{fig:SSTPredictions}
\end{figure}

\begin{figure}
    \centering
    \subfloat{
    \begin{overpic}[width=0.93\linewidth]{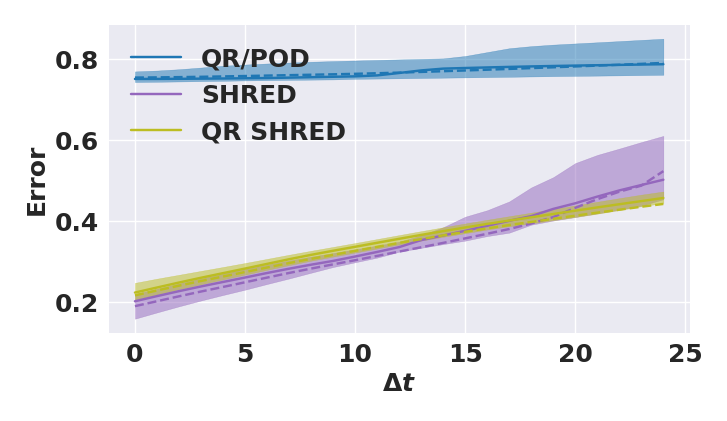}
      \put(2,55){\textbf{A}}
    \end{overpic}}
    
    \subfloat{
    \begin{overpic}[width=0.93\linewidth]{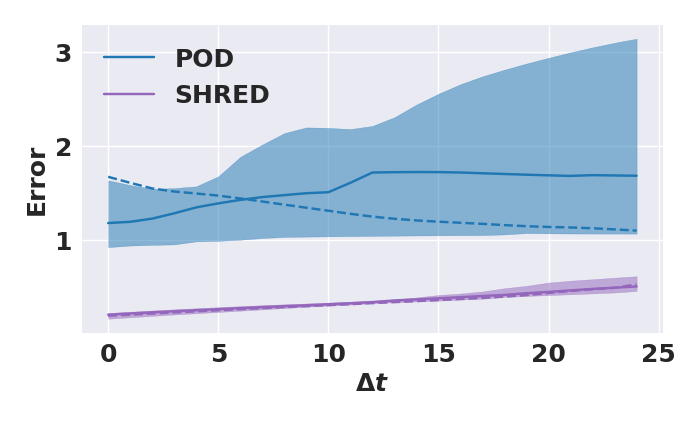}
      \put(2,55){\textbf{B}}
    \end{overpic}}
    
    \caption{Forecasting error for turbulent flow with $\Delta t = 0.006$s. A LSTM is trained to forecast QR placed, panel \textit{A}, or randomly placed, panel \textit{B}, sensor measurements which are subsequently used to perform reconstructions using SHRED and gappy POD. In this experiment, the validation set is selected from earlier in the data so that the forecast occurs immediately after following training data. 16 forecasts are performed and the median error is denoted by the solid lines, with the 25th and 75th percentiles of reconstruction error defining the shaded region.  The dashed line represents an ensembled forecast.}
    \label{fig:ISOPredictions}
\end{figure}
\subsection{Forecasting of high-dimensional spatio-temporal fields}

\label{sec:results2}

In this section, we explore how the expressiveness of SHRED models with few sensors can be leveraged to perform forecasts of high-dimensional spatio-temporal states using only a limited subsampling of the state.  Forecasting the evolution of high-dimensional states from sensor measurements is an exceedingly challenging task, owing to the fact that doing so combines the problem of system identification with the difficulties of forecasting.  We propose a two step approach that leverages the success of LSTMs for low-dimensional time-series forecasting \cite{lindemann_survey_2021} and SHRED for state estimation.

Let $\{ \bx _i \}_{1}^{T_t}$ represent a training dataset of temporally ordered states with corresponding sensor measurements $\{ \by _i \} _{1}^{T_t}$.  We train an LSTM network, $\mathcal{G} \left( \{ \by_i \}_{t- k}^{t}\right)$, to map between a trajectory of sensor measurements and the subsequent measurement, $\by_{t+1}$, by finding

\begin{equation}
    \mathcal{G} \in \argmin_{\widetilde{\mathcal{G}} \in  \mathscr{G}} \sum _{i=1}^N ||\by_{t+1} - \widetilde{\mathcal{G}}\left( \{ \by _j \}_{i-k}^i \right)||_2
\end{equation}

using the ADAM optimizer on the training data.  Doing so, we then forecast beyond the time interval of the training data to obtain forecasted sensor measurements $\{ \hat{\by}_i\}_{T_t + 1}^{T_t + p}$ for $p > 1.$  The forecasted measurements are used by a SHRED model, trained on $\{ \bx _i \}_{1}^{T_t}$, to obtain forecasts of the high-dimensional state $\bx.$ Detailed training protocols are included in SI. 

The forecasted states $\{\hat{\bx}_i \}_{T_t +1} ^{T_t + p}$ are evaluated against the ground truth for each $\Delta t$ forecast rather than the mean error across all forecasts to demonstrate performance over forecasts of varying length. The error at time $p\Delta t$ is thus given by
\begin{equation}
    \text{Error} = \frac{ \left| \left|  \mathcal{H} \left( \{\hat{\by} _j \}_{p-k}^p\right) - \bx_p \right| \right|_2}{|| \bx _p ||_2}. 
\end{equation}
We also include reconstructions from forecasted sensor measurements obtained by gappy POD for comparison, akin to the work of \cite{parish_time-series_2020}.  Finally, we only show results for SST and turbulent flow data, as the compressed form of the atmospheric ozone data is biased towards both QR placement and POD reconstructions.

For the SST example, we select the first 85\% of the dataset to act as training data, the subsequent 20 snapshots as validation data, and the remainder as the test set.  A LSTM for forecasting is trained on the training data and sensor measurements are forecast beyond the validation set.  These forecasted sensor measurements are then used to construct a forecast of the high-dimensional spatio-temporal field using a trained SHRED model and gappy POD.  We consider the cases that the sensors are placed via QR or randomly.  We perform 16 runs for each experiment, and consider an ensembled forecast found by averaging the forecast of each run.  The results for QR placed and randomly placed sensors are shown in Fig. \ref{fig:SSTPredictions} with $\Delta t$ representing one week.  With QR placement, the forecasts obtained by SHRED outperform that of POD in the short-term and are similar for longer forecast horizons.  However, with randomly placed sensors SHRED is still able to obtain comparably accurate forecasts while POD fails to consistently yield faithful reconstructions.

Unlike sea-surface temperature, even medium-range forecasts of a turbulent flow are impossible due to the fact that the flow is not quasi-periodic or stationary in nature.  For this reason, we focus on achieving accurate short-term forecasts in this application.  To do so, we use the same scheme as before, with the exception that the validation set is selected to occur earlier in the training data.  Of the 1567 snapshots of the turbulent flow with sufficiently long preceding measurement histories, the first 1000 are selected as training, followed by 50 samples selected for validation.  The next 100 samples are also used for training and the remainder constitute the test set.  This scheme allows for more accurate short-term forecasts because the forecast occurs directly after the training data, as opposed to directly after the validation set.  Fig. \ref{fig:ISOPredictions} shows the results for the cases that sensors are placed via QR and randomly with $\Delta t$ representing 0.006 seconds. In both cases, the forecast obtained by SHRED deteriorates quickly but greatly outperforms that obtained by POD reconstructions.

\subsection{Temporal spacing of time-delays}
\begin{figure*}[t]
    \centering
    \includegraphics[width=0.95\textwidth]{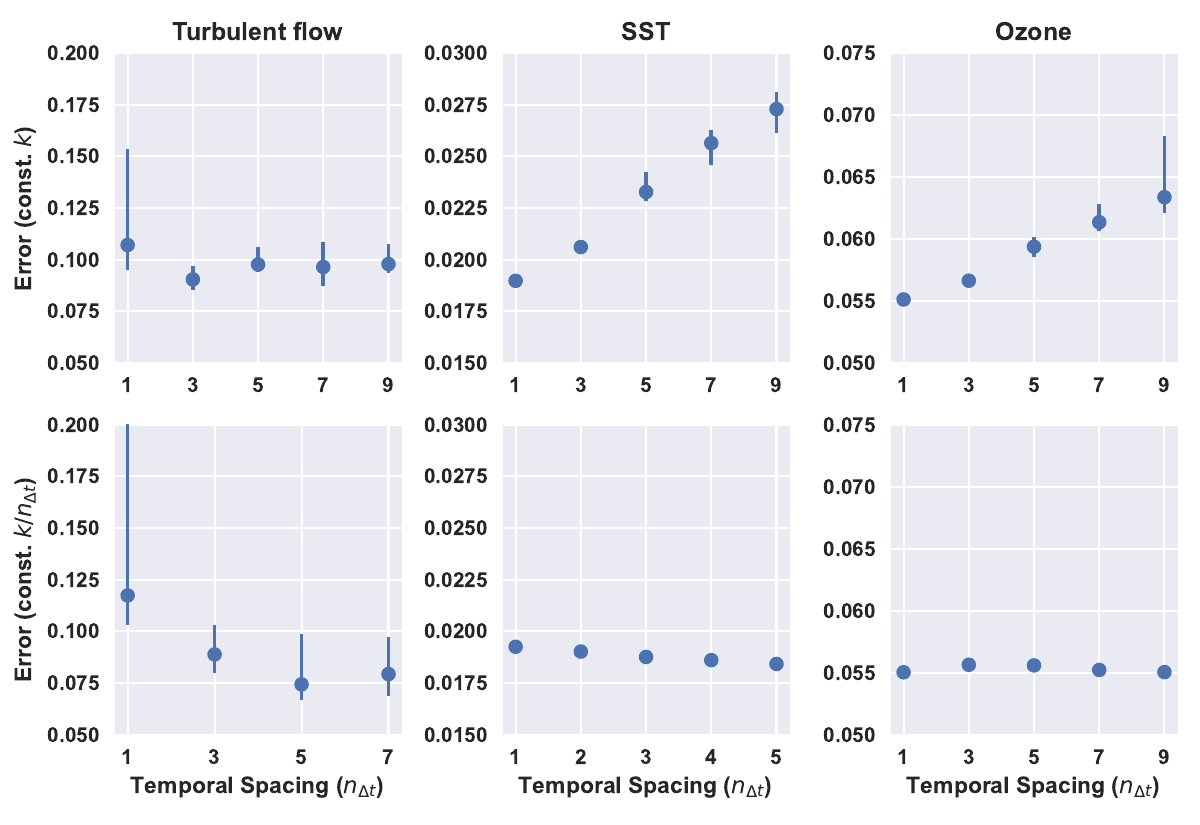}
    \caption{Error results for SHRED models in the three considered datasets when the temporal spacing of the trajectories used for reconstruction is varied. $k$ denotes the maximal time-delay  of the trajectories used by a SHRED model. $n_{\Delta t}$ denotes the temporal spacing between utilized time-delays in terms of of the discretization of the data, $\Delta t$. The top row considers the case where $k$ is approximately fixed to the values specified in Section \ref{sec:results1} and the bottom row considers the case where $k$ varies in addition to $n_{\Delta t}$ such that $k/ n_{\Delta t}$ is fixed to the values of $k$ specified in  \ref{sec:results1}. Each plot shows the median and interquartile range of reconstruction error from 32 trained SHRED models.} 
    \label{fig:time-delays}
\end{figure*}

    To this point, we have only considered SHRED models where the temporal spacing between utilized sensor measurements is equal to the discretization, $\Delta t$, of the data snapshots. However, time-delay embedding techniques can often be sensitive to the selection of both the maximal time-delay incorporated, $k\Delta t$, and the temporal spacing of those time-delays, $n_{\Delta t} \Delta t$ \cite{abarbanel_analysis_2012}. In this notation, Sections \ref{sec:results1} and \ref{sec:results2} show results for $n_{\Delta t} = 1$ and in this section we consider the performance of SHRED for varying $n_{\Delta t}.$

    We perform two sets of experiments for each considered dataset. The first varies $n_{\Delta t}$ but keeps $k$ approximately constant. This corresponds to a decrease in the total number of sensor measurements in the input to a SHRED model. The second set of experiments varies $n_{\Delta t}$ and $k$ such that the total number of utilized sensor measurements is held constant. As a result, the first $N - k n_{\Delta t}$ snapshots in each dataset must be discarded. We therefore modify the number of samples in the training, validation, and test sets such that there are 900 training samples for the turbulent flow, 1000 for the SST, and 1000 for the atmospheric ozone data. The results for both experiments are reported in Fig. \ref{fig:time-delays}. In both cases, all other network hyperparameters and the manner in which the train/validation/test split is performed remain unchanged from Section \ref{sec:results1}.

    We find that the performance of SHRED in the cases of SST and atmospheric ozone concentration remains largely unchanged when $n\Delta t$ is varied, provided $k$ is also scaled to use the same number of sensor measurements. If $n_{\Delta t}$ is increased and $k$ is held fixed, performance deteriorates for larger values of $n_{\Delta t}.$

    On the other hand, our results are less clear when SHRED is applied to the turbulent flow data. It appears that commensurately increasing $n_{\Delta t}$ and $k$ moderately improves performance, while increasing only $n_{\Delta t}$ does not have a strong effect on reconstruction accuracy. However, it is difficult to draw strong conclusions from these results as there are strong indications that 900 training samples are insufficient to train the model. This is illustrated by moderately worse median performance and much larger variance for the cases $n_{\Delta t} = 1$ in Fig. \ref{fig:time-delays}. The only difference between these models and those from Section \ref{sec:results1} is a reduction in access to training data. The effect of this reduction is likely more visible in the turbulent flow data than SST or atmospheric ozone because it is intrinsically a more complex dataset.

\section{Discussion}

We have demonstrated SHRED as a method for system identification/state estimation and the forecasting of high-dimensional time-series data. By including sensor trajectories in the model for reconstruction, SHRED outperformed competing techniques based on static reconstructions from gappy POD or shallow decoder networks. We considered three high-dimensional, spatio-temporal datasets, a synthetically generated forced turbulent flow, a simulation of atmospheric ozone concentration from a chemical transport model, and weekly mean sea-surface temperature.  

The superior performance of SHRED held across both interpolatory (reconstruction) and extrapolatory (forecasting) regimes. In Section \ref{sec:results1}, the task was purely an interpolatory reconstruction. Training and test data, upon which we evaluate reconstruction accuracy, were drawn from the same distribution.  As a result, overfitting to the training data was not an issue, and SHRED was able to capture fine grain details of complex flows with as few as one or three sensors. Applied in this manner, SHRED offers an attractive method for compression of large datasets as well as a framework for developing reduced order models of dynamical systems. The performance of SHRED in this section is also indicative of the promise of SHRED for high-dimensional data that is stationary or periodic. Notably, SHRED demonstrated remarkable indifference towards sensor placement. In this interpolatory regime, issues arising from the necessity of optimal sensor placement can be mitigated by increasing the information content of inputs to a neural network through the use of sensor trajectories.  These results show that the sensor trajectories encode a significant amount of information, just as is expected of time-delayed embeddings of dynamics~\cite{brunton2017chaos,arbabi2017ergodic,lee2019linking,pan2020structure,kamb2020time,bakarji2022discovering}.

In Section \ref{sec:results2}, the task at hand was forecasting the evolution of high-dimensional spatio-temporal data.  In contrast to Section \ref{sec:results1}, training and test sets are temporally separated and thus the results are extrapolatory in nature. High-dimensional forecasts were accomplished by forecasting low-dimensional sensor measurements and constructing states from the predicted measurements. Forecasting is an inherently difficult task for complex systems, especially in the high-dimensional fields often encountered in the engineering and the physical sciences. Still, SHRED’s use of sensor trajectories allows for improved forecasting performance as compared to existing methods. Furthermore, SHRED for forecasting exhibits robustness to sensor placement not found in the other methods considered. We thus advocate for the use of SHRED for both the tasks of state estimation and state prediction from forecasted sensor measurements. Further work must be done to improve the accuracy of low-dimensional sensor forecasts to be used in conjunction with SHRED for state prediction. Improving the forecasts of high-dimensional data could enable on-the-fly compression and novel sensing mechanisms.

Finally, future work will explore the extension of SHRED to multi-modal data, as briefly discussed in the section on coupled PDEs. This extension can be viewed as an augmentation of the full-state representation to include sensor measurements from various modalities. The inputs to the network then have the interpretation of encoding both the spatial location and sensing modality. As in the cases considered here, the application of SHRED in multi-modal settings will require a high-fidelity model with priors that well-represent the physical system of interest. Provided such a model and coupling between the field of observations and augmented full-state, in principle it is possible to apply our methodology to reconstruct the full-state from measurements of one of its constituent fields.

\section*{Data Availability}
Each of the three datasets considered in this paper have been released in previous work \cite{li_public_2008, reynolds_improved_2002, velegar_scalable_2019}. For details of accessing the data, we direct the reader to the SI.

\section*{Code Availability}
The code used to generate the results in this paper can be found at (https://github.com/Jan-Williams/pyshred).

\section*{Acknowledgements}

The authors acknowledge support from the National Science Foundation AI Institute in Dynamic Systems (grant number 2112085).  JNK further acknowledges support from the Air Force Office of Scientific Research (FA9550-19-1-0011 and FA9550-19-1-0386).

\bibliographystyle{naturemag}
\bibliography{mainv5.bib}

\end{document}

%% file: z_vardef.tex
\newcommand{\bb}{\mathbf{b}}
\newcommand{\bc}{\mathbf{c}}

\newcommand{\bv}{\mathbf{v}}

\newcommand{\bx}{\mathbf{x}}

\newcommand{\by}{\mathbf{y}}
\newcommand{\bh}{\mathbf{h}}

\newcommand{\bC}{\mathbf{C}}

\newcommand{\bQ}{\mathbf{Q}}
\newcommand{\bR}{\mathbf{R}}

\newcommand{\bU}{\mathbf{U}}
\newcommand{\bV}{\mathbf{V}}
\newcommand{\bW}{\mathbf{W}}
\newcommand{\bX}{\mathbf{X}}

\DeclareMathOperator*{\argmax}{arg\rm{}max}
\DeclareMathOperator*{\argmin}{arg\rm{}min}